\documentclass{ws-fnl}
\usepackage{graphicx,url,amsmath,amssymb,pgfornament} 

\title{Prolegomena to the Bestiary}
\author{Yang-Hui He}
\address{London Institute for Mathematical Sciences, Royal Institution, W1S 4BS, UK \\
    Merton College, University of Oxford, OX14JD, UK\\
hey@maths.ox.ac.uk}

\newcommand{\ornamentSep}{\noindent\hfil{{\pgfornament[width=0.618033988749895\columnwidth,color=black]{88}}}}

\begin{document}

\maketitle

\begin{abstract}
``Calabi-Yau Manifolds: a Bestiary for Physicists'' by Tristan Hubsch  in 1992 was a classic that served to introduce algebraic geometry to physicists when the first string theory revolution of 1984 - 94 brought, inter alia, the subject of Calabi-Yau manifolds to the staple of high-energy theorists.
We are fortunate that a substantially expanded and updated new edition of the Bestiary will shortly appear.
This brief note will serve as an afterword to the much anticipated volume.
\end{abstract}

~\\
~\\

When Professor Tristan H\"ubsch and World Scientific asked me to write a few words for the publication of the Second Edition of {\it Calabi-Yau Manifolds: A Bestiary For Physicists}, I felt tremendously honoured.
It was as if an Aquinian monk had been summoned to scribe a commentative addendum to the {\it Summa}, or a neo-Confucian scholar, to the {\it Analects}.
In my generation of mathematical physicists - pursuing PhDs at the turn of the millennium - the {\it Bestiary} was a unique and valued book.
If anything, it was perhaps the only textbook on our shelves - where the typical title would be ``Introduction to $x$'', where $x$ is some technical word in theoretical physics or pure mathematics - to have so exotic a word as ``bestiary''. Even the font (Gothic) and the cover picture (photo of origami and notebook) added an eccentric but inviting charm.

The book, as it turned out, was exactly that; it was {\it inviting}.
In the early decades of the twentieth century, the ground-breaking work of Einstein and Hilbert made it clear that physicists should learn differential geometry.
By the mid-decades, that of Dirac, Wigner, Yang, et al.~did likewise for abstract algebra and representation theory.
In the last decades, that of Atiyah, Chern, Witten, et al.- highlighted perhaps by Witten's Fields Medal in 1990 - inspired theorists that algebraic geometry, too, should be an essential tool in a physicist's machinery.
But, there's the rub.

To quote a perspicacious comment from Terrence Tao in a recent interview \cite{tao}, ``\ldots in the sixties and seventies \ldots there was an era of mathematics where abstraction was making huge strides in simplifying and unifying a lot of mathematics that was previously very empirical \ldots this is sometimes what we call the Bourbaki era in mathematics. And it did veer a little bit too far from being grounded.''
This Bourbaki style resulted in tomes of the revered GTM yellow \cite{gtm}, definitive but intimidating, rigorous but impenetrable, serving far more as decorative references on the shelves to impress a visitor than to initiate a neophyte into the subject.
In algebraic geometry, at the time of the 1980-90s when physicists realized the necessity to study this amazing field, there was Hartshorne's canon \cite{hartshorne}, which was considered even amongst graduate students of pure mathematics as an Olympian challenge. For physicists, it might as well have been written in Sumerian.
The friendlier alternatives would have been Reid \cite{reid}, Harris \cite{harris}, or \cite{abhyankar} (which did not have the advanced material requisite to gauge theory), or Griffith-Harris \cite{gh} (832 pages!).
Thus, as wonderful as these books are, they were not what was desperately needed by the theoretical physics community. 

H\"ubsch's {\it Bestiary} was the saving grace.
Though the title suggested a rather technical focus on Calabi-Yau manifolds, it really was an invitation to complex algebraic geometry for the non-specialist.
As esoteric as the word ``Calabi-Yau'' might sound, there is important historical context.
In the nascence of topology, Euler classified (smooth, compact, orientable, boundary-less) surfaces by a single integer, the {\it genus} $g$, which characterizes the number of ``holes''.
Thereafter, Gauss, in his Theorema Egregium, related the Euler number $\chi = 2 - 2g$ to the integral of curvature, whereby linking topology to analysis.
Here, we see a {\bf trichotomy}: (i) $g=0$ for $\chi>0$ (positive curvature); (ii) $g=1$ for $\chi=0$ (zero curvature); and (iii) $g>1$ for $\chi<0$ (negative curvature).
For surfaces, there are respectively the sphere, the torus, and those of higher genera (including the famous anti-de-Sitter space in real dimension two, as well as modular curves as quotients of the upper-half complex plane).
The intermediate case of zero-curvature obviously distinguishes itself. Indeed, the torus is precisely a (the) Calabi-Yau manifold of complex dimension 1.

Generalizing this story to beyond surfaces had been a central problem in mathematics since Euler-Gauss-Riemann. At least in complex geometry (even-dimensional manifolds that admit a complex structure), and in particular those with an additional K\"ahler structure, the above trichotomy persists.
Through the index theorem, the Euler number $\chi$ is given by the integral of the Chern class (of the tangent bundle), giving us 
(i) Fano varieties of  positive curvature; (ii) Calabi-Yau varieties of zero curvature; and (iii) varieties of general type.
The uniqueness and conjectural existence of the K\"ahler metric that determines the Chern class, by Calabi \cite{calabi} in the 1950s, and the subsequent Fields-winning demonstration of the existence by Yau \cite{yau} in the late 1970s, were great triumphs in modern mathematics.

All beautiful mathematics must find their place in physics.
In the 1980s, when string theorists came across the conditions of ``Ricci-flat'' and ``K\"ahler'' in their search for a consistent, supersymmetric, low-energy effective theory from the 10-dimensional heterotic string \cite{chsw,hull} in order to attain the Standard Model of particle physics, Strominger's office was next to Yau's at the IAS. Geometry then unknown to physicists was at once available from the expert who just obtained the Fields Medal. The two worlds coincided, and sparkled ideas and collaborations that would persist till the very present. Indeed, the name ``Calabi-Yau'' was coined by the physicists in honour of the two great mathematicians.

While Yau gave some immediate examples of Calabi-Yau manifolds on which one might compactify string theory, even meeting some of the rudimentary constraints such as getting the right number of generations of fermions turned out to be a non-trivial challenge.
How to find (smooth, compact) Calabi-Yau 3-folds (complex dimension 3, real dimension 6) that satisfied physical constraints which are translated into exact conditions in algebraic and differential geometry?
The hunt was on.

From the earliest cases of complete intersections in a single projective space, to a single hypersurface \cite{wp4} in weighted projective $\mathbb{C}\mathbb{P}^4$, to complete intersections in products in projective spaces -- the so-called CICYs \cite{cicy,cicy2,cicy3,Green:1987cr,Kim:1989dq}, the number of explicitly constructed Calabi-Yau manifolds rapidly rose to the thousands.
The fact that it was the physicists who had to pursue such constructions proved very much the aforementioned point about the predominant air of abstraction in the pure mathematics community left from the Bourbaki era.
Somehow it was considered - a sentiment which is luckily no longer the main-stream today - uncivilized to stoop from abstraction to concrete examples in mathematics.
In a way, the {\it Bestiary} documented this flurry of activity in the last decade of the twentieth century and threw down the gauntlet back at the Bourbakians.

\ornamentSep
~\\

Some three decades have passed since the publication of the first edition of the {\it Bestiary} \cite{bestiary1}.
Enormous strides, of course, have since been made in physics and in mathematics.
The {\it First String Revolution}, wherein the different formulations of string theory were shown to be dual to each other, with the compactifications of which leading to Calabi-Yau manifolds, have been roughly marked by historians to be circa 1984-94.
This was the inspiration of the {\it Bestiary}.
Its last chapter alluded to {\bf mirror symmetry}, where the complex and K\"ahler structures of a given Calabi-Yau manifold $M$ should be exchanged.
While in physics, this followed from a duality between super-conformal field theories \cite{mirror}, the implication to geometry, where complex and symplectic structures could be mixed, was astounding.
Mirror symmetry is now a vast field in itself and the reader is referred to \cite{mirror1,mirror2,mirror3,mirror4} for more recent textbooks.

The search for mirror pairs was a parallel skein in mathematics to the search for the Standard Model from Calabi-Yau compactifications in physics.
The search led to the remarkable works of Batyrev-Borisov in complete intersections in toric varieties coming from reflexive polytopes \cite{BB}, the tour-de-force classification thereof \cite{Kreuzer:1995cd}, and the more recent works on generalized complete intersections \cite{Anderson:2015iia,GWinf}.
All at once, there is a simple combinatorial criterion for mirror symmetry for large classes of manifolds, and the number of explicit smooth Calabi-Yau 3-folds grew to billions.

The {\it second string revolution} came in the period of 1994 - 2003, when M-theory, D-branes and the AdS/CFT holographic correspondence enriched strings to extended objects.
For Calabi-Yau manifolds, it opened up vast interests in non-compact and singular ones (of course, M-theory also ignited interests in $G_2$ manifolds, many of which can be constructed from fibrations over Calabi-Yau 3-folds).
These ``local'' Calabi-Yau manifolds are key to AdS/CFT \cite{Maldacena:1997re} because the D-branes can ``probe'' them \cite{Klebanov:1998hh,Hanany:1998sd}.
In this geometric context, the profound holographic principle relating quantum gravity in dimension $d$ to gauge theory on its $(d-1)$-dimensional boundary exhibits as the correspondence of the algebraic geometry of a moduli space and the representation theory of quivers with potential.
The reader is referred to pedagogical books in \cite{holo1,holo2,holo3} for the physics, and to \cite{Bao:2020sqg} for some of the mathematics, especially to the McKay Correspondence \cite{mckayADE} and to bi-partite structures \cite{Yamazaki:2008bt,bipartite}.

The {\it third string revolution} started in 2003 when the {\bf string landscape} and the statistics of vacuum configurations were the concepts to seep into the minds of scientists and the public alike.
Prompted by the plethora of the Kreuzer-Skarke Calabi-Yau manifolds, as well as the increasing number of further structures one could place on these manifolds in compactification - such as stable vector bundles, D-branes wrapping algebraic cycles, etc. - Doulgas \cite{stats} and KKLT \cite{Kachru:2003aw} contemplated, on a concrete footing, the vastness of possible four-dimensional universes obtainable from string theory.

Whilst astronomical numbers such as $10^{500}$ began to dishearten the public confidence in string theory, it must be pointed out that such numbers came from back-of-the-envelop estimates (typical Betti numbers of Calabi-Yau manifolds in a generic calculation). More importantly, the whole point of the landscape (and related swampland \cite{Ooguri:2006in}) is not to plague us with a plethora, but 
to rule out the overwhelming majority of scenarios.
Indeed, even without imposing cosmological constraints and the correct Yukawa couplings, reduction via the constraint of getting the right particle content alone should be by some 8 orders of magnitude \cite{HetMSSM,Gmeiner:2005vz,Candelas:2007ac,hetLine,Candelas:2016fdy}.
Thus, though there seems to be an unfathomably large landscape of possible four-dimensional universes obtainable from string theory, so far, a single one that is consistent with our observations - the Standard Model of particles, general relativity dictating large-scale universe with a small positive cosmological constant (de-Sitter space), and quantum gravity governing black-holes - is yet to be found.

In mathematics, this situation is curiously paralleled.
Recall the trichotomy discussed above, of positive (Fano), zero (Calabi-Yau), and negative curvature (general type) complex manifolds.
The Fano case is known to have a finite number of topological types in every dimension.
The general type is known to be wildly unclassifiable and infinite in number.
For the boundary case of Calabi-Yau, it is a conjecture of Yau that in every dimension the number of different topological types (Betti numbers, Hodge numbers, intersection numbers, etc.) is finite. In complex dimension one, we see this is true: there is only the torus $T^2$.
In complex dimension two, there are only the 4-torus $T^4$ and the K3-surface.
In complex dimension three (and above), despite the vast numbers, one still suspects finite possibilities in topology.
The statement is even stronger \cite{reidfantasy}: the moduli space of (smooth, compact, simply-connected) Calabi-Yau 3-folds is fantasized to be connected, so that all of them are related by bi-rational transformations such as flops.

~\\
We see, therefore, Calabi-Yau manifolds still very much occupy a central place in contemporary mathematics and theoretical physics (even completely outside of string theory, Calabi-Yau structures are unexpectedly manifesting in the Standard Model \cite{He:2014oha} and in ordinary Feynman amplitudes \cite{Bourjaily:2018ycu}).
Irrespective of one's interests in Calabi-Yau spaces in any context, it should be noted that in many ways, the {\it Bestiary} was not merely a pioneer in inviting researchers to algebraic geometry, but to {\it computational} algebraic geometry.
This brings us to the final, and by far, the most important aspect of our discussion.

The second half of the twentieth century has convinced us of the indispensable {\it r\^ole of computers} in mathematics (broadly defined to include any theoretical research) \cite{icm2018}.
The speeding up of computations and simulations certainly go without saying, but more pertinent are the ever-increasing cases where major fundamental results {\it could not} be attained without the computer. The four-colour theorem, the Kepler Conjecture, the classification theorem of finite simple groups, to name but a few, could not have been possible to humans alone, simply because the proofs or case-checking would take longer than a life-span to go through.

Indeed, the data of Calabi-Yau manifolds were compiled since the 1990s by impressive computer work. For instance, the CICYs, arguably the first database ever created for geometry, were found on the CERN super-computer. 
With the astounding growth and availability of the personal computer, already by the Kreuzer-Skarke list in 2000, the half-billion manifolds were found on a Pentium desktop.
Since then, there is a host of databases and freely-available software dedicated to studying the Calabi-Yau (and related) landscape, such as:
\begin{description}
    \item[The Graded Ring Database ]
    \url{http://www.grdb.co.uk/}
    \item[Kreuzer-Skarke ] 
    \url{http://hep.itp.tuwien.ac.at/~kreuzer/CY/}
    \item[KS (some triangulations/orientifolds) ]
    \url{http://www.rossealtman.com/toriccy/}
    \item[KS Explorer ]
    \url{https://cy.tools/}
    \item[Picard-Fuchs operators ]
    \url{https://cydb.mathematik.uni-mainz.de/}
    \item[CICY+ (at bottom)]
    \url{https://www.physics.ox.ac.uk/our-people/lukas}
    \item[Various resources, e.g., numerical metrics, etc. ] ~\\
    \url{https://github.com/topics/calabi-yau-manifolds}
\end{description}

Meanwhile, the paucity of friendly introductions to algebraic geometry at the time of \cite{bestiary1} is well supplanted by a number of excellent textbooks that emphasize on explicit calculations rather than abstruse sophistry \cite{CLO,m2book,toricbook}, as well as some dedicated to physicists and engineers \cite{bookComp,appliedGeo,abhyankar}.

~\\

Something was in the air in 2017, and four independent groups thought of using machine-learning to study problems in string compactifications \cite{He:2017aed,Krefl:2017yox,Ruehle:2017mzq,Carifio:2017bov}.
In particular, \cite{He:2017aed} proposed that algebraic geometry, if not any mathematical problem, could be phrased as an image-processing query for a neural network and that the structure of mathematics aught to be studied by AI \cite{He:2021oav}.
A flurry of activity on AI-guided theoretical discovery in related and unrelated areas in the past lustrum is an encouraging testimony that Calabi-Yau manifolds reside at the confluence of pure mathematics, theoretical physics, and AI research.
One illustrative time-line might give the reader an informative glimpse. There have been a number of branch-offs of the annual string theory conference, which since 1986 has been bringing together various sub-communities in high-energy theory and in mathematics.
The following are the beginning dates of these international conference series:
\begin{description}
 \item[1986-]
  {\it ``Strings''}
 
 \item[2002-]
  {\it ``StringPheno''}: string theorists interested in obtaining phenomenologically viable physics

\item[2006 - 2010]
{\it String Vacuum Project}: an NSF funded multi-institution collaboration for the string landscape
 
 \item[2008 -]
 {\it ``ISGT''} Integrability in String/Gauge Theory
 
 \item[2011-]
  {\it ``String-Math''}: interactions between string theory and mathematics

 \item[2012-]
  {\it ``Amplitudes''}: computing QFT amplitudes using string-inspired techniques

 \item[2014-] 
  {\it SIAM annual meeting} with a String/Theoretical Physics session

 \item[2017-]
  {\it ``String-Data''}: machine-learning for the string landscape

 \item[2020-]
  {\it ''DANGER''}: (Data in ANalysis, GEometry \& Representation Theory) 
\end{description}

The trend toward data and AI driven research is self-evident.
We refer the reader to the first textbook \cite{He:2018jtw} and editorial \cite{MLbook} on machine-learning in theoretical physics and pure mathematics, as well as a wonderful pedagogical monograph on modern data scientific methods \cite{Ruehle:2020jrk} for physicists (cf.~\cite{He:2023csq,Berglund:2023ztk,Hirst:2023kdl}).

In many ways, \cite{He:2018jtw}, focusing on machine-learning for the Calabi--Yau landscape, as an attempt to invite three communities (mathematics, physics, computer-science/AI) to a cocktail party, was my homage to the {\it Bestiary}, which had been so influential to my intellectual youth.
I am thrilled to see the appearance of a new edition, with extensive additions and revisions by Prof. H\"ubsch, and am flattered to be given the opportunity to write these words.
I certainly hope the current generation of theoretical physicists and pure mathematicians would embrace the {\it Bestiary} with the same excitement and curiosity that I did as a student.

\section*{Acknowledgements}
I am indebted to STFC for grant ST/J00037X/2, and the Leverhulme Trust for a project grant.
Indeed, I am grateful to the London Institute and Merton College, Oxford, for providing two corners of paradise, wherein I could contemplate without distraction.

\newpage

\end{document}